\documentclass{amsart}[12pt]
\usepackage{amsmath, amsthm, amscd, amsfonts}

\begin{document}

\title{A REVIEW ON SOME GEOMETRIC RESULTS OF THE \u{S}MULIAN'S THEOREM ON FRECHET DIFFERENTIABILITY OF NORMS }

\maketitle

\author \hspace{1.15cm} A. ASSADI,H. HAGHSHENAS AND H. HOSSEINI GUIVE

\begin{center}
Department of Mathematics, Birjand University,Iran\\
E-mail:h$_{-}$haghshenas60@yahoo.com
\end{center}

\begin{abstract}
In this paper, we prove the \u{S}mulian's theorem on Fr\'{e}chet
differentiability of norm, [3], and present some of its geometric
results concerning the Gateaux and Fr\'{e}chet differentiability
of norm and properties of the allied space and its dual such as
reflexivity and strict convexity.
\end{abstract}AMS Subject Classification(2000) : 46B20\\
\textbf{\textit{Key Words:}}\hspace{0.5cm} strictly convex space,
locally uniformly convex space, Gateaux and Fr\'{e}chet
differentiability.
\section{Introduction}Differentiability of convex continuous
functions has been studied extensively during the recent fifty
years. Among all functions in this class, the norm function is of
the great importance and deserves to be studied in more detail,
since the geometric structure of the space is highly related to
it.In this paper we emphasize on Differentiability of norm and use
the \u{S}mulian's theorem to prove some geometric results.\\
\section{Basic Definitions and Preliminaries}In this section we
recall some elementary definitions and theorems required for the
next parts of the text. Through this paper, $X$ is a normed vector
space with norm $\|.\|$ and the unit sphere $S(X)$.\\\\
\textbf{Definition 1. }$X$ is said to be \textit{locally uniformly
convex} if for each sequence $(x_{n})_{n=1}^{\infty}$ and any
point $x$ in $S(X)$, $\displaystyle{\lim_{n\rightarrow
\infty}\|x_{n}+x\|=2}$ implies $\displaystyle{\lim_{n\rightarrow
\infty}\|x_{n}-x\|=0}$. Moreover, $X$ is \textit{strictly convex}
if $x=y$ whenever $x,y\in S(X)$ and $(x+y)/2\in S(X)$.\\\\
\textbf{Remark 1. }Clearly, each locally uniformly convex space is
strictly convex.\\For example, the Euclidean space $\Re^{2}$ is a
locally uniformly convex space.\\\\\textbf{Definition 2.   }The
norm function, $\|.\|$, is Gateaux differentiable at $0\neq x$ if
there is some $x^{*}\in X^{*}$ such that the quantity
$\displaystyle{\lim_{t\rightarrow 0}\frac{\|x+ty\|-\|x\|}{t}}$
exists for each $y\in X$ and is equal to $\langle x^{*},y
\rangle$. If the limit exists uniformly for each $y\in S(X)$, then
the norm function is said to be Fr\'{e}chet differentiable at
$x$.\\\\The following theorems are proved in
[2]:\\\\\textbf{Theorem 1.   } if $X^{*}$ is strictly convex then
the norm function $\|.\|$ on $X$ is Gateaux differentiable
.\\\\\textbf{Theorem 2. }The norm function is Fr\'{e}chet
differentiable at $0\neq x\in X$ if and only if
$$\displaystyle{\lim_{t\rightarrow 0}\frac{\|x+ty\|+\|x-ty\|-2\|x\|}{t}=0,}$$
uniformly for each $y\in S(X)$.\section{Main Results}We are now
ready to prove the \u{S}mulian's theorem and state some of its
geometric results.\\\\\textbf{Theorem 3. }(\u{S}mulian). Suppose
$x\in S(X)$. The following are equivalent:\\\textbf{a. }The norm
function is Fr\'{e}chet differentiable at $x$.\\\textbf{b.   }For
all $(f_{n})_{n=1}^{\infty},(g_{n})_{n=1}^{\infty}\subseteq
S(X^{*})$, if $\displaystyle{\lim_{n\rightarrow
\infty}f_{n}(x)=1}$ and $\displaystyle{\lim_{n\rightarrow
\infty}g_{n}(x)=1}$, then $\displaystyle{\lim_{n\rightarrow
\infty}\|f_{n}-g_{n}\|=0}$.\\\textbf{c. }Each
$(f_{n})_{n=1}^{\infty}\subseteq S(X^{*})$ with
$\displaystyle{\lim_{n\rightarrow \infty}f_{n}(x)=1}$ is
convergent in $S(X^{*})$.\\\\\textit{Proof.
   }$\textbf{a}\Rightarrow \textbf{b})$ Let the norm function is Fr\'{e}chet differentiable
   at $x$. For each $\varepsilon >0$ there exists $\delta >0$ such
   that $$\|x+y\|+\|x-y\|\leq 2+\varepsilon\|y\|,$$ for each $y$
   with $\|y\|<\delta$. Suppose $(f_{n})_{n=1}^{\infty},(g_{n})_{n=1}^{\infty}\subseteq
   S(X^{*})$ with $\displaystyle{\lim_{n\rightarrow \infty}f_{n}(x)= \displaystyle{\lim_{n\rightarrow
   \infty}g_{n}(x)=1}}$. Since $\displaystyle{\lim_{n\rightarrow
   \infty}f_{n}(x)=1}$, there exists $N_{1}\in \aleph$ such that
   $|f_{n}(x)-1|<\varepsilon\delta$ for each $n\geq N_{1} $.
   Similarly, there is $N_{2}\in \aleph$ such that $|g_{n}(x)-1|<\varepsilon\delta$ for each $n\geq N_{2}
   $. Let $M = max\{ N_{1},N_{2}\}$. Then for any $n>M$,
   $$|f_{n}(x)-1|<\varepsilon\delta \ \ \ and \ \ \ |g_{n}(x)-1|<\varepsilon\delta$$
   and we have:
 \begin{eqnarray*}
 (f_{n}-g_{n})(y)&=&f_{n}(x+y)+g_{n}(x-y)-f_{n}(x)-g_{n}(x)\\
 &\leq& \|x+y\|+\|x-y\|-f_{n}(x)-g_{n}(x)\\
 &\leq& 2+\varepsilon\|y\|-f_{n}(x)-g_{n}(x)\\
 &\leq&|f_{n}(x)-1|+|g_{n}(x)-1|+\varepsilon\|y\|\\
 &\leq&3\varepsilon\delta.
 \end{eqnarray*}
   Hence for each $n>M$,
   $$\|f_{n}-g_{n}\|= \displaystyle{\sup_{h\in S(X)}(f_{n}-g_{n})(h)= \displaystyle{\sup_{h\in S(X)}
   \frac{(f_{n}-g_{n})(\delta h)}{\delta}\leq 3\varepsilon.}}$$ $\textbf{b}\Rightarrow \textbf{a})$
   Suppose norm is not Fr\'{e}chet differentiable at $x$. So there
   is $\varepsilon >0$ such that for each $\delta=\frac{1}{n}$, $\|x+y_{n}\|+\|x-y_{n}\|>2+\varepsilon\|y_{n}\|$
   for some $y_{n}$ with $\|y_{n}\|\leq\frac{1}{n}$. Now according
   to the Hahn-Banach extension theorem, for each $n\in \aleph$,
   there are $(f_{n})_{n=1}^{\infty},(g_{n})_{n=1}^{\infty}\subseteq
   S(X^{*})$ such that
   $$f_{n}(x+y_{n})=\|x+y_{n}\|\hspace{1cm} and \hspace{1cm}g_{n}(x-y_{n})=\|x-y_{n}\|.$$
   Also for each $n\in
   \aleph$,$$|\|x+y_{n}\|-\|x\||\leq\|y_{n}\|\leq\frac{1}{n}.$$Thus
   $|f_{n}(x+y_{n})-1|<\frac{1}{n}$ which implies that
   $\displaystyle{\lim_{n\rightarrow\infty}f_{n}(x+y_{n})=1}$. On the other
   hand, $|f_{n}(y_{n})|\leq\|y_{n}\|\leq\frac{1}{n}$ which
   follows that $\displaystyle{\lim_{n\rightarrow\infty}f_{n}(y_{n})=0}$.
   Therefore $\displaystyle{\lim_{n\rightarrow\infty}f_{n}(x)=1}$. In a similar
   way, $\displaystyle{\lim_{n\rightarrow\infty}g_{n}(x)=1}$.Since $$g_{n}(x)\leq\|g_{n}\|=1
   \hspace{1cm}and\hspace{1cm}f_{n}(x)\leq\|f_{n}\|=1,$$ we have
   $f_{n}(x)+g_{n}(x)\leq2$. Thus we conclude that
   $$(f_{n}-g_{n})(y_{n})\geq\|x+y_{n}\|+\|x-y_{n}\|-2>\varepsilon\|y_{n}\|.$$
   Consequently $\|f_{n}-g_{n}\|\geq\varepsilon$ which is a
   contradiction.\\\\$\textbf{b}\Rightarrow \textbf{c})$ Let $(f_{n})_{n=1}^{\infty}\subseteq
   S(X^{*})$ such that $\displaystyle{\lim_{n\rightarrow\infty}f_{n}(x)=1}$.
   There is some $f\in S(X^{*})$ such that $f(x)=\|x\|=1$. For
   each $n\in \aleph$, let $g_{n}=f$. Then
   $\displaystyle{\lim_{n\rightarrow\infty}g_{n}(x)=f(x)=1}$. Hence using \textbf{b} we
   find that $\displaystyle{\lim_{n\rightarrow\infty}\|f_{n}-g_{n}\|=0}$ which
   means that $(f_{n})_{n=1}^{\infty}$ is convergent.\\\\$\textbf{c}\Rightarrow \textbf{b})$
    Consider $(f_{n})_{n=1}^{\infty},(g_{n})_{n=1}^{\infty}\subseteq
    S(X^{*})$ with
    $$\displaystyle{\lim_{n\rightarrow\infty}f_{n}(x)= \displaystyle{\lim_{n\rightarrow\infty}g_{n}(x)=1}}.$$For
    each $y\in X$, define $h_{n}(y)=f_{\frac{n+1}{2}}(y)$ for odd values of $n$ and
    $h_{n}(y)=g_{\frac{n}{2}}(y)$for even values of $n$. So we
    have $$h_{2n}=g_{n}\in S(X^{*})\hspace{1cm}and\hspace{1cm}h_{2n-1}=f_{n}\in
    S(X^{*}),$$ and therefore $h_{n}\in S(X^{*})$. Moreover,
    $$h_{2n}(x)=g_{n}(x)\rightarrow1\hspace{1cm}and\hspace{1cm}h_{2n-1}(x)=f_{n}(x)\rightarrow1.$$
    Hence $h_{n}(x)\rightarrow1$. Using \textbf{c} there is $h\in
    S(X^{*})$ such that $\|h_{n}-h\|\rightarrow0$. Then,
    $$\|f_{n}-g_{n}\|=\|h_{2n-1}-h_{2n}\|\leq
    \|h_{2n-1}-h\|+\|h_{2n}-h\|\rightarrow0$$,
    as required.\\\\\textbf{Corollary 1.   }If $X^{*}$ is locally
    uniformly convex, then the norm function is Fr\'{e}chet
    differentiable on $X$.\\ \textit{Proof.   }Suppose $x\in
    S(X)$, $f\in S(X^{*})$, $f(x)=1$ and $(f_{n})_{n=1}^{\infty}\subseteq
    S(X^{*})$ such that $\displaystyle{\lim_{n\rightarrow\infty}f_{n}(x)=1}$. Then
    $2\geq\|f_{n}+f\|\geq(f_{n}+f)(x)\rightarrow2 $ as $n\rightarrow\infty$.
    Therefore$$\displaystyle{\lim_{n\rightarrow\infty}(2\|f_{n}\|^{2}+2\|f\|^{2}-\|f_{n}+f\|^{2})=0.}$$
    Now Since $X^{*}$ is locally uniformly convex $\displaystyle{\lim_{n\rightarrow\infty}\|f_{n}-f\|=0}$. Thus
    \u{S}mulian's theorem completes the proof.\\\\\textbf{Corollary 2.
    }If $X^{*}$ is Fr\'{e}chet differentiable, then $X$ is
    reflexive.\\ \textit{Proof.   }It is known that $X$ is
    reflexive if and only if each nonzero $f\in X^{*}$ attains its
    norm at some $x\in S(X)$. Let $f\in S(X^{*})$ and choose $(x_{n})_{n=1}^{\infty}\in S(X)$
    such that $f(x_{n})\rightarrow1$. By the \u{S}mulian's
    theorem, $\displaystyle{\lim_{n\rightarrow\infty}x_{n}=x\in S(X)}$. Therefore,
    $$f(x)=f(\displaystyle{\lim_{n\rightarrow\infty}x_{n})= \displaystyle{\lim_{n\rightarrow\infty}f(x_{n})=1=\|f\|.}}$$
    If now $f\in X^{*}$ is non-zero, then $\frac{f}{\|f\|}\in
    S(X^{*})$ and according to the above manner there exists $x\in
    S(X)$ such that $\frac{f}{\|f\|}(x)=1$.\\\\ \textbf{Corollary 3.   } Let $X$
    is finite dimensional. If the norm is Gateaux differentiable,
    then it is Fr\'{e}chet differentiable.\\ \textit{Proof.
    }Since the norm is Gateaux differentiable at $x\in S(X)$,
    there exists a unique $f\in S(X^{*})$ such that $f(x)=1$. Let $(f_{n})_{n=1}^{\infty}\subseteq
   S(X^{*})$ and $\displaystyle{\lim_{n\rightarrow\infty}f_{n}(x)=1}$. According
   to the \u{S}mulian's theorem, it is sufficient to show
   $(f_{n})_{n=1}^{\infty}$ converges. Since $S(X^{*})$ is a closed
   and bounded subset in a finite dimensional space it is compact.
   Thus $(f_{n})_{n=1}^{\infty}$ has a convergent subsequence
   $(f_{n_{k}})_{k=1}^{\infty}$. Let
   $\displaystyle{\lim_{k\rightarrow\infty} f_{n_{k}}=f_{1}}$. Consequently,
    $$\displaystyle{\lim_{k\rightarrow\infty} f_{n_{k}}(x)=f_{1}(x)=1.}$$But
    there is only one $f\in S(X^{*})$ such that $f(x)=1$. So
    $f=f_{1}$. Therefore each $(f_{n})_{n=1}^{\infty}\subseteq
   S(X^{*})$ has a subsequence convergent to $f$ which follows
   that $(f_{n})_{n=1}^{\infty}$ converges to $f$.\\

\end{document}